\documentclass [10pt] {article}

 \usepackage [cp1251] {inputenc}
 \usepackage [english] {babel}
 \usepackage {amssymb}
 \usepackage {mathtools}
 \usepackage {amsthm}
 \usepackage [curve] {xypic}
 \newdir {:=}{{}}      

 \let \"\~             
 \let \acute\'
 \def \={\bar}         
 \def \~{\tilde}
 \def \^{\hat}
 \def \smash{\wedge}   
 \def \rdto{\mathbin|} 
 \def \divs{\mathrel|} 
 \def \le{\leqslant}
 \def \ge{\geqslant}
 \def \:{\colon}       
 \def \qto{\dashrightarrow}        
 \def \<{\langle}      
 \def \>{\rangle}
 \def \[{\lvert}       
 \def \]{\rvert}
 \def \({\lVert}       
 \def \){\rVert}
 \let \$\#             
 \def \+{\prescript{+}{}}          
 \def \`{\lfloor}
 \def \'{\rfloor}
 \newcommand* {\ga} [2] {#1\lfloor#2\rfloor}
 \newcommand* {\ai} [2] {{\rceil}#1\lfloor#2\rfloor}
 \def \%{{\rceil}}     
 \renewcommand* {\#} [1] {\lvert#1\rvert}      
 \renewcommand* {\*} [1] {\lvert#1\rvert}      
 \newcommand* {\bc} [2]            
              {\Big(\negthickspace
               \begin{array}{c}#1\\#2\end{array}
               \negthickspace\Big)}
 \def \|{\prescript}   

 \def \N{{\boldsymbol{\mathsf N}}}
 \def \Z{{\boldsymbol{\mathsf Z}}}
 \def \Q{{\boldsymbol{\mathsf Q}}}

 \def \R{{\mathcal R}} 

 \def \K{{\mathcal K}}   

 \def \E{{\mathcal E}} 

 \def \id{{\mathrm{id}}}

 \DeclareMathOperator {\Deg} {Deg}
 \DeclareMathOperator {\Hom} {Hom}
 \DeclareMathOperator {\im} {im}
 \DeclareMathOperator {\Fix} {Fix}
 \DeclareMathOperator {\sAr} {sAr}
 \DeclareMathOperator {\qAr} {\tilde sAr}

 \DeclareMathOperator {\proj} {pr} 
 \DeclareMathOperator {\incl} {in} 
 \def \M{M}

 \newcommand* {\head} [1]
 {\subsubsection* {\mathversion{bold}#1}}

 \newcommand* {\subhead} [1]
 {\addvspace\medskipamount
  \noindent {\mathversion{bold}\bf\itshape #1\/}}

 \newenvironment* {claim} [1] []
 {\begin{trivlist}\item [\hskip\labelsep {\bf #1}] \it}
 {\end{trivlist} }

 \newenvironment* {demo} [1] []
 {\begin{trivlist}\item [\hskip\labelsep {\it #1}] }
 {\end{trivlist} }

\begin {document}

 \title {\large\bf
         On homotopy invariants of finite degree}

 \author {\normalsize\rm
          S.~S.~Podkorytov}

 \date {}

 \maketitle

 \begin {abstract} \noindent
 We prove that
 homotopy invariants of finite degree
 distinguish homotopy classes of maps
 of a connected compact CW-complex
 to a nilpotent connected CW-complex
 with finitely generated homotopy groups.
 \end {abstract}


 \head {\S~1. Introduction}


 $\N=\{0,1,\dotsc\}$.
 {\it Space\/} means ``pointed topological space''.
 {\it CW-complexes\/} are also pointed
 (the basepoint being a vertex).
 {\it Map\/} means ``basepoint preserving continuous map''.
 Homotopies,
 the notation $[X,Y]$, etc.\
 are to be understood in the pointed sense.


 \subhead {Invariants of finite degree.}
 Let
 $X$ and $Y$ be spaces,
 $V$ be an abelian group, and
 $f\:[X,Y]\to V$ be a function (a homotopy invariant).
 Let us define a number $\Deg f\in\N\cup\{\infty\}$,
 the {\it degree\/} of $f$.
 Given a map $a\:X\to Y$ and a number $r\in\N$,
 we have the map $a^r\:X^r\to Y^r$ (the Cartesian power),
 which induces the homomorphism
 $C_0(a^r)\:C_0(X^r)\to C_0(Y^r)$
 between the groups of (unreduced) zero-dimensional chains
 with the coefficients in $\Z$.
 Let the inequality $\Deg f\le r$ be equivalent to
 the existence of a homomorphism
 $l\:\Hom(C_0(X^r),C_0(Y^r))\to V$ such that
 $f([a])=l(C_0(a^r))$ for all maps $a\:X\to Y$.
 As one easily sees,
 $\Deg f$ is well defined by this condition.
 {\it Finite-degree\/} invariants are those of finite degree.


 \subhead {Main results.}

 \begin {claim} [1.1. Theorem.]
 Let
 $X$ be a connected compact CW-complex,
 $Y$ be a nilpotent connected CW-complex
 with finitely generated homotopy groups, and
 $u_1,u_2\in[X,Y]$ be distinct classes.
 Then,
 for some prime $p$,
 there exists a finite-degree invariant $f\:[X,Y]\to\Z_p$
 such that
 $f(u_1)\ne f(u_2)$.
 \end {claim}

 Related facts were known for certain cases
 where $[X,Y]$ is an abelian group \cite{me-sp, me-st}.
 Theorem~1.1 follows
 (see \S~11)
 from a result of Bousfield--Kan and Theorem~1.2.

 We call a group {\it $p$-finite\/}
 (for a prime $p$)
 if
 it is finite and
 its order is a power of $p$.

 \begin {claim} [1.2. Theorem.]
 Let
 $p$ be a prime,
 $X$ be a compact CW-complex, and
 $Y$ be a connected CW-complex
 with $p$-finite homotopy groups.
 Then every invariant $f\:[X,Y]\to\Z_p$ has finite degree.
 \end {claim}

 Probably,
 Theorem~1.2 can be deduced from Shipley's convergence theorem
 \cite{Shi},
 which we do not use.
 We use an (approximate) simplicial model of $Y$
 that admits a harmonic (see \S~6) embedding in a simplicial
 $\Z_p$-module.


 \subhead {Non-nilpotent examples.}
 The following examples show the importance of the nilpotency
 assumption in Theorem~1.1.
 We consider finite-degree invariants on $\pi_n(Y)=[S^n,Y]$.

 \begin {claim} [1.3.]
 Let $Y$ be a space with $\pi_1(Y)$ perfect.
 Then,
 for any abelian group $V$,
 any finite-degree invariant $f\:\pi_1(Y)\to V$ is constant.
 \end {claim}

 \begin {demo}
 This follows from Lemmas 12.2 and 3.6.
 \qed
 \end {demo}

 \begin {claim} [1.4.]
 Take $n>1$.
 Let $Y$ be a space such that
 $\pi_n(Y)\cong\Z^2$ and
 an element $g\in\pi_1(Y)$ induces an order 6 automorphism on
 $\pi_n(Y)$.
 Then,
 for any abelian group $V$,
 any finite-degree invariant $f\:\pi_n(Y)\to V$ is constant.
 \end {claim}

 \begin {demo}
 This follows from Lemmas 12.2 and 12.3 and claim~3.7.
 \qed
 \end {demo}


 \subhead {An example: maps $S^{n-1}\times S^n\to S^n_{(\Q)}$
           \rm (cf.\ \cite[Example~4.6]{ArkLup}).}
 Take an even $n>0$.
 Let
 $c\:S^{n-1}\times S^n\to S^{2n-1}$ be a map of degree 1.
 Put
 $i=[\id]\in\pi_n(S^n)$,
 $j=i*i\in\pi_{2n-1}(S^n)$ (the Whitehead square), and
 $u(q)=(qj)\circ[c]\in[S^{n-1}\times S^n,S^n]$, $q\in\Z$.
 Let $l\:S^n\to S^n_\Q$ be the rationalization.
 Put $\=u(q)=[l]\circ u(q)\in[S^{n-1}\times S^n,S^n_\Q]$.
 The classes $u(q)$, $q\in\Z$, are pairwise distinct;
 moreover,
 the classes $\=u(q)$, $q\in\Z$, are pairwise distinct
 (the proof is omitted).

 Is it true that,
 under the assumptions of Theorem~1.1,
 there must exist an $r\in\N$ such that
 the elements of $[X,Y]$ are distinguished by invariants of
 degree at most $r$?
 No,
 as the following claim shows.

 \begin {claim} [1.5.]
 Let
 $V$ be an abelian group and
 $f\:[S^{n-1}\times S^n,S^n]\to V$ be an invariant of degree at
 most $r\in\N$.
 Then $f(u(q))=f(u(0))$ whenever $r!\divs q$.
 \end {claim}

 The following claim shows the importance of the assumption of
 Theorem~1.1 that $Y$ has finitely generated homotopy groups.

 \begin {claim} [1.6.]
 Let
 $V$ be an abelian group and
 $f\:[S^{n-1}\times S^n,S^n_\Q]\to V$ be an invariant of finite
 degree.
 Then $f(\=u(q))=f(\=u(0))$, $q\in\Z$.
 \end {claim}

 The following claim shows that,
 under the assumptions of Theorem~1.1,
 finite-degree invariants taking values in $\Q$
 may not distinguish rationally distinct homotopy classes.

 \begin {claim} [1.7.]
 Lat $f\:[S^{n-1}\times S^n,S^n]\to\Q$ be an invariant of
 finite degree.
 Then $f(u(q))=f(u(0))$, $q\in\Z$.
 \end {claim}


 \subhead {Elusive elements of $H_0(Y^X)$.}
 The space of maps $X\to Y$ is denoted $Y^X$.
 An invariant $f\:[X,Y]\to V$ gives rise to the homomorphism
 $\+f\:H_0(Y^X)\to V$, $\`u\'\mapsto f(u)$
 (here $\`u\'$ denotes the basic element corresponding to $u$).
 Is it true that,
 under the assumptions of Theorem~1.1,
 for any non-zero element $w\in H_0(Y^X)$
 there exist
 an abelian group $V$ and
 a finite-degree invariant $f\:[X,Y]\to V$
 such that
 $\+f(w)\ne0$?
 No,
 as the following claim shows.

 \begin {claim} [1.8.]
 Take $n>1$.
 Let
 $Y$ be a space and
 $u_1,u_2\in\pi_n(Y)$ be elements of coprime finite orders.
 Put $w=\`u_1+u_2\'-\`u_1\'-\`u_2\'+\`0\'$.
 Let
 $V$ be an abelian group and
 $f\:\pi_n(Y)\to V$ be an invariant of finite degree.
 Then $\+f(w)=0$.
 \end {claim}

 \begin {demo}
 This follows from Lemmas 12.2 and 3.8.
 \qed
 \end {demo}

 If the group $\pi_n(Y)$ is torsion and divisible,
 then the same is true for any elements $u_1,u_2\in\pi_n(Y)$
 (this follows from Lemmas 12.2 and 3.9).
 In this case,
 $\pi_n(Y)$ cannot be finitely generated
 (without being zero).
 In return,
 $Y$ can be $p$-local,
 e.~g.\
 $Y=\K(P,n)$ (the Eilenberg--MacLane space)
 for $P=\Z[1/p]/\Z$.


 \head {\S~2. Preliminaries}


 We say
 {\it crew\/} for ``pointed set'' and
 {\it archism\/} for ``basepoint preserving function''.
 We use the standard model structure on the category of
 simplicial crews (and archisms) \cite[Corollary~3.6.6]{Hov}.
 The words {\it fibration}, {\it cofibration}, etc. refer to
 it.
 A {\it fibring\/} simplicial archism is a fibration.
 An {\it isotypical\/} simplicial archism,
 or an {\it isotypy},
 is a weak equivalence.
 {\it Isotypic\/} simplicial crews are weakly equivalent ones.

 An abelian group is a crew
 (the basepoint being $0$);
 a simplicial abelian group is a simplicial crew.

 We call a simplicial crew $T$
 {\it compact\/} if it is generated by a finite number of
 simplices, and
 {\it gradual\/} if the crews $T_q$, $q\in\N$, are finite.

 For simplicial crews $K$ and $T$,
 we have a simplicial crew $T^K$, the function object
 (denoted $\hom_*(K,T)$ in \cite[Ch.~VIII, 4.8]{BouKan}).
 A simplicial archism $f\:K\to L$ induces
 a simplicial archism $T^f\:T^L\to T^K$, etc.
 We use this notation in the topological case as well.

 The sign $\sim$ denotes the homotopy relation;
 the sign $\simeq$ denotes the homotopy eqivalence of spaces.


 \subhead {Main homomorphisms.}
 By default,
 chains and homology have coefficients in a commutative ring
 $\R$;
 $\Hom=\Hom_\R$.
 (In \S~1, we had $\R=\Z$ implicitly.)

 For spaces $X$ and $Y$,
 define $\R$-homomorphisms
 $$
 \|XY\mu_r\:C_0(Y^X)\to\Hom(C_0(X^r),C_0(Y^r)),
 \qquad
 \`a\'\mapsto C_0(a^r),
 $$
 $r\in\N$.
 We have the projection $\R$-homomorphism
 $$
 \|XY\nu\:C_0(Y^X)\to H_0(Y^X).
 $$

 For simplicial crews $K$ and $T$,
 define $\R$-homomorphisms
 $$
 \|KT\mu_r\:C_0(T^K)\to\Hom_0(C_*(K^r),C_*(T^r)),
 \qquad
 \`b\'\mapsto C_*(b^r),
 $$
 $r\in\N$.
 Here
 $\`b\'$ is the basic chain corresponding to a simplex
 $b\in(T^K)_0$,
 i.~e.\ a simplicial archism $b\:K\to T$;
 $b^r\:K^r\to T^r$ is the Cartesian power;
 $C_*(b^r)\:C_*(K^r)\to C_*(T^r)$ is the induced
 $\R$-homomorphism of graded $\R$-modules of chains;
 $\Hom_0$ denotes the $\R$-module of grading-preserving
 $\R$-homomorphisms.
 We have the projection $\R$-homomorphism
 $$
 \|KT\nu\:C_0(T^K)\to H_0(T^K).
 $$


 \head {\S~3. Group algebras and gentle functions}


 Let $\ga\R G$ denote the group $\R$-algebra of a group $G$.
 An element $g\in G$ has the corresponding basic element
 $\`g\'\in\ga\R G$.
 The augmentation ideal $\ai\R G\subseteq\ga\R G$ is the kernel
 of the $\R$-homomorphism $\ga\R G\to\R$, $\`g\'\mapsto1$.
 The ideal $\ai\R G^s$ ($s>0$) is $\R$-generated by elements
 of the form $(1-\`g_1\')\dotso(1-\`g_s\')$.

 Let $V$ be an abelian group.
 A function $f\:G\to V$ gives rise to the homomorphism
 $\+f\:\ga\Z G\to V$, $\`g\'\mapsto f(g)$.
 We call $f$
 {\it $r$-gentle\/}
 if $\+f\rdto\ai\Z G^{r+1}=0$, and
 {\it gentle\/} (or {\it polynomial\/})
 if it is $r$-gentle for some $r\in\N$
 \cite[Ch.~V]{Pas}.

 Let $p$ be a prime.

 \begin {claim} [3.1. Lemma.]
 Let $U$ be a finite $\Z_p$-module of dimension $m$.
 Then $\ai{\Z_p}U^{(p-1)m+1}=0$.
 \qed
 \end {claim}

 \begin {claim} [3.2. Corollary.]
 Let $U$ and $V$ be $\Z_p$-modules.
 If $U$ is finite,
 then every function $f\:U\to V$ is gentle.
 \qed
 \end {claim}

 \begin {claim} [3.3. Lemma \rm {\cite[Proposition~1.2]{Dre}}.]
 Let
 $U$, $V$, and $W$ be abelian groups,
 $f\:U\to V$ be an $r$-gentle function, and
 $g\:V\to W$ be an $s$-gentle one 
 ($r,s\in\N$).
 Then the function $g\circ f\:U\to W$ is $rs$-gentle.
 \end {claim}

 \begin {demo}
 This follows from \cite[Ch.~V, Theorem~2.1]{Pas}.
 \qed
 \end {demo}

 A function $f\:U\to V$ between abelian groups induces the
 $\R$-homomorphism $f_\R\:\ga\R U\to\ga\R V$,
 $\`u\'\mapsto\`f(u)\'$.

 \begin {claim} [3.4. Corollary.]
 Let
 $U$ and $V$ be abelian groups and
 $f\:U\to V$ be an $r$-gentle ($r\in\N$) function.
 Then,
 for any $s\in\N$,
 the $\R$-homomorphism $f_\R$ maps
 the ideal $\ai\R U^{rs+1}$ to the ideal $\ai\R V^{s+1}$.
 \qed
 \end {claim}

 \begin {claim} [3.5. Lemma.]
 Let $I$ be a set.
 For each $i\in I$,
 let
 $U_i$ and $V_i$ be abelian groups and
 $f_i\:U_i\to V_i$ be an $r$-gentle ($r\in\N$) function.
 The the function
 $$
 \prod_{i\in I} f_i\:
 \prod_{i\in I} U_i\to
 \prod_{i\in I} V_i
 $$
 is $r$-gentle.
 \qed
 \end {claim}


 The following claims are used
 only in discussion of the examples of \S~1,
 not in the proof of the main results.

 \begin {claim} [3.6. Lemma.]
 Let
 $G$ be a perfect group and
 $V$ be an abelian group.
 Then any gentle function $f\:G\to V$ is constant.
 \end {claim}

 \begin {demo}
 This follows from \cite[Ch.~III, Corollary~1.3]{Pas}.
 \qed
 \end {demo}

 \begin {claim} [3.7.]
 Let
 $U$ be an abelian group isomorphic to $\Z^2$,
 $J\:U\to U$ be an automorphism of order 6,
 $V$ be an abelian group, and
 $f\:U\to V$ be a gentle function.
 Suppose that
 the function $\Z\times U\to V$, $(t,u)\mapsto f(J^tu-u)$, is
 gentle.
 Then $f$ is constant.
 \end {claim}

 \begin {demo}
 The proof is omitted.
 \qed
 \end {demo}

 \begin {claim} [3.8. Lemma.]
 Let
 $U$ and $V$ be abelian groups,
 $f\:U\to V$ be a gentle function, and
 $u_1,u_2\in U$ be elements of coprime finite orders.
 Then $f(u_1+u_2)-f(u_1)-f(u_2)+f(0)=0$.
 \qed
 \end {claim}

 \begin {claim} [3.9. Lemma.]
 Let
 $U$ be a divisible torsion abelian group, and
 $V$ be an abelian group.
 Then every gentle function $f\:U\to V$ is 1-gentle.
 \qed
 \end {claim}

 \begin {claim} [3.10. Lemma.]
 Let $G$ and $H$ be groups.
 Then the ideal $\ai\R{G\times H}^s$ ($s>1$) is $\R$-generated
 by elements of the form $(1-\`a_1\')\dotso(1-\`a_{s-q}\')
 (1-\`b_1\')\dotso(1-\`b_q\')$,
 where
 $0\le q\le s$,
 $a_t\in G\times1\subseteq G\times H$, and
 $b_t\in 1\times H\subseteq G\times H$.
 \qed
 \end {claim}

 \begin {claim} [3.11. Lemma.]
 A function $F\:\Z\to\Q$ is $r$-gentle ($r\in\N$)
 if and only if
 it is given by a polynomial of degree at most $r$.
 \qed
 \end {claim}


 \head {\S~4. Keys of a commutative square}


 Let $E$ be a commutative ring.
 Consider the diagram of simplicial $E$-modules and
 $E$-homomorphisms
 $$
 \xymatrix {
 V'
 \ar[d]^-{f'}
 \ar@/^2ex/[rr]^-{t'} &&
 W
 \ar[ll]^-{g'}
 \ar[d]_-{g''} \\
 U
 \ar@/^2ex/[u]^-{s'}
 \ar@/_2ex/[rr]_-{s''} &&
 V'',
 \ar[ll]_-{f''}
 \ar@/_2ex/[u]_-{t''}
 }
 $$
 where the square is commutative: $f'\circ g'=f''\circ g''$.
 We call the quadruple $(s',s'',t',t'')$ a {\it key\/} of this
 square if
 we have $(-s',s'')\circ(-f',f'')+(g',g'')\circ(t',t'')=\id$ in
 the diagram
 $$
 \xymatrix {
 U
 \ar@/_2ex/[rr]_-{(-s',s'')} &&
 V'\oplus V''
 \ar[ll]_-{(-f',f'')}
 \ar@/_2ex/[rr]_-{(t',t'')} &&
 W.
 \ar[ll]_-{(g',g'')}
 }
 $$
 The pair $(t',t'')$ is called a {\it half-key\/} in this case.

 \begin {claim} [4.1. Lemma.]
 Let
 $$
 \xymatrix {
 V'
 \ar[d]^-{f'}
 \ar@/^2ex/[rr]^-{t'} &&
 W
 \ar[ll]^-{g'}
 \ar[d]_-{g''} \\
 U &&
 V''
 \ar[ll]_-{f''}
 \ar@/_2ex/[u]_-{t''}
 }
 $$
 be a commutative square of simplicial $E$-modules and
 $E$-homomorphisms with a half-key,
 $T$ be a simplicial crew, and
 $k'\:T\to V'$ and $k''\:T\to V''$ be simplicial archisms such
 that $f'\circ k'=f''\circ k''$.
 Consider the simplicial archism
 $l=t'\circ k'+t''\circ k''\:T\to W$.
 Then
 $g'\circ l=k'$ and
 $g''\circ l=k''$.
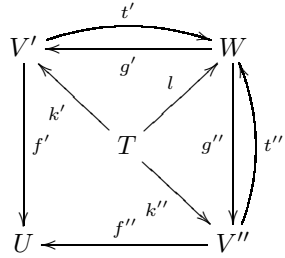
 $$
 \xymatrix {
 V'
 \ar[dd]^-{f'}
 \ar@/^2ex/[rr]^-{t'} &&
 W
 \ar[ll]^-{g'}
 \ar[dd]_-{g''} \\
 &
 T
 \ar[ul]^-{k'}
 \ar[dr]_-{k''}
 \ar[ur]^-{l} &
 \\
 U &&
 V''
 \ar[ll]_-{f''}
 \ar@/_2ex/[uu]_-{t''}
 }
 $$
 \qed
 \end {claim}

 By a {\it sector\/} of a simplicial $E$-homomorphism
 $h\:\~W\to W$ we mean a simplicial $E$-homomorphism
 $s\:W\to\~W$ such that $h\circ s=\id$.

 \begin {claim} [4.2. Lemma.]
 Consider a commutative diagram of simplicial $E$-modules
 and $E$-homomorphisms
 $$
 \xymatrix {
 0 &&
 \~U
 \ar[ll]
 \ar[d]^-{f} &&
 \~V
 \ar[ll]_-{\~p}
 \ar[d]^-{g} &&
 \~W
 \ar[ll]_-{\~q}
 \ar[d]^-{h} &&
 0
 \ar[ll] \\
 0 &&
 U
 \ar[ll] &&
 V
 \ar[ll]_-{p} &&
 W
 \ar[ll]_-{q} &&
 0.
 \ar[ll]
 }
 $$
 Suppose that
 its rows are split exact and
 $h$ has a sector.
 Then the left-hand square has a key.
 \end {claim}

 \begin {demo} [Proof.]
 Let $(k,l)$ and $(\~k,\~l)$
 (see the diagram below)
 be splittings:
 \begin {align*}
 p\circ k & =\id, &
 l\circ q & =\id, &
 k\circ p+q\circ l & =\id, \\
 \~p\circ\~k & =\id, &
 \~l\circ\~q & =\id, &
 \~k\circ\~p+\~q\circ\~l & =\id,
 \end {align*}
 and
 $s$ be a sector: $h\circ s=\id$.
 Put
 $r=\~q\circ s\circ l$ and
 $\^k=\~k+r\circ(k\circ f-g\circ\~k)$.
 Then $(0,k,\^k,r)$ is a key.
 $$
 \xymatrix {
 0 &&
 \~U
 \ar[ll]
 \ar[d]^-{f}
 \ar@/_2ex/[rr]_-{\~k}
 \ar@/^3ex/[rr]^-{\^k} &&
 \~V
 \ar[ll]_-{\~p}
 \ar[d]_-{g}
 \ar@/_2ex/[rr]_-{\~l} &&
 \~W
 \ar[ll]_-{\~q}
 \ar[d]_-{h} &&
 0
 \ar[ll] \\
 0 &&
 U
 \ar[ll]
 \ar@/_2ex/[rr]_-{k}
 \ar@/^2ex/[u]^-{0} &&
 V
 \ar[ll]_-{p}
 \ar@/_2ex/[rr]_-{l}
 \ar@/_2ex/[u]_-{r} &&
 W
 \ar[ll]_-{q}
 \ar@/_2ex/[u]_-{s} &&
 0
 \ar[ll]
 }
 $$
 \qed
 \end {demo}


 \begin {claim} [4.3. Lemma.]
 Let
 $L$ and $M$ be simplicial crews,
 $j\:L\to M$ be an isotypical cofibration, and
 $Q$ be a fibrant simplicial crew.
 Then $Q^j\:Q^M\to Q^L$ is an isotypical fibration.
 \qed
 \end {claim}

 \begin {claim} [4.4. Lemma.]
 Let
 $Q$ and $R$ be simplicial crews,
 $c\:Q\to R$ be a fibration, and
 $N$ be a simplicial crew isotypic to a point.
 Then $c^N\:Q^N\to R^N$ is an isotypical fibration.
 \qed
 \end {claim}

 \begin {claim} [4.5. Lemma.]
 Suppose that $E$ is a field.
 Let
 $V$ and $W$ be simplicial $E$-modules and
 $f\:W\to V$ be an isotypical fibring simplicial
 $E$-homomorphism.
 Then $f$ has a sector.
 \qed
 \end {claim}

 \begin {claim} [4.6. Lemma.]
 Suppose that $E$ is a field.
 Let
 $L$ and $M$ be simplicial crews,
 $j\:L\to M$ be an isotypical cofibration,
 $Q$ and $R$ be simplicial $E$-modules, and
 $c\:Q\to R$ be a fibring simplicial $E$-homomorphism.
 Then the commutative square
 $$
 \xymatrix {
 Q^L
 \ar[d]_-{c^L} &&
 Q^M
 \ar[ll]_-{Q^j}
 \ar[d]^-{c^M} \\
 R^L &&
 R^M
 \ar[ll]_-{R^j}
 }
 $$
 has a key.
 \end {claim}

 \begin {demo} [Proof.]
 Consider the (strictly) cofibration sequence
 $$
 \xymatrix {
 L
 \ar[rr]^-{j} &&
 M
 \ar[rr]^-{k} &&
 N.
 }
 $$
 Since $j$ is isotypical,
 the simplicial crew $N$ is isotypic to a point.
 We have the following diagram of simplicial $E$-modules and
 $E$-homomorphisms:
 $$
 \xymatrix {
 0 &&
 Q^L
 \ar[ll]
 \ar[d]^-{c^L} &&
 Q^M
 \ar[ll]_-{Q^j}
 \ar[d]^-{c^M} &&
 Q^N
 \ar[ll]_-{Q^k}
 \ar[d]^-{c^N} &&
 0
 \ar[ll] \\
 0 &&
 R^L
 \ar[ll] &&
 R^M
 \ar[ll]_-{R^j} &&
 R^N
 \ar[ll]_-{R^k} &&
 0.
 \ar[ll]
 }
 $$
 We show that the rows are split exact.
 Consider the upper row.
 Obviously,
 it is exact in the middle and the right-hand terms.
 $Q$ is fibrant
 since it is a simplicial abelian group.
 By Lemma~4.3,
 $Q^j$ is an isotypical fibration.
 By Lemma~4.5,
 $Q^j$ has a sector.
 Therefore,
 the upper row is split exact.
 The same is true for the lower row.
 By Lemma~4.4,
 $c^N$ is an isotypical fibration.
 By Lemma~4.5,
 $c^N$ has a sector.
 By Lemma~4.2,
 the desired key exists.
 \qed
 \end {demo}


 \head {\S~5. Quasisimplicial archisms}


 A {\it quasisimplicial archism\/} $f\:K\qto L$
 between simplicial crews $K$ and $L$
 is a sequence of archisms $f_q\:K_q\to L_q$, $q\in\N$.
 Let
 $\qAr(K,L)$ denote the crew of quasisimplicial archisms and
 $\sAr(K,L)$ denote the subcrew of simplicial ones.

 A quasisimplicial archism $f\:U\qto V$
 between simplicial abelian groups is {\it $r$-gentle\/} if
 the archisms $f_q\:U_q\to V_q$ are $r$-gentle.

 Let $T$ be a simplicial crew.
 For $m,q\in\N$,
 let
 $[m|q]$ be the set of non-strictly increasing functions
 $[m]\to[q]$
 (where $[q]=\{0,\dotsc,q\}$) and
 consider the archism
 $$
 T(m,q)=(T(h))_{h\in[m|q]}\:T_q\to T_m^{[m|q]}.
 $$
 We call $T$ {\it $m$-soluble\/} if,
 for any $q$,
 the archism $T(m,q)$ is injective.

 Let $p$ be a prime.

 \begin {claim} [5.1. Lemma.]
 Let
 $T$ be a gradual simplicial crew,
 $U$ be a gradual simplicial $\Z_p$-module,
 $R$ be an $m$-soluble ($m\in\N$) simplicial $\Z_p$-module,
 $d\:T\to U$ be a cofibration, and
 $k\:T\to R$ be a simplicial archism.
 Then,
 for some $r\in\N$,
 there exists an $r$-gentle quasisimplicial archism
 $w\:U\qto R$ such that
 $w\circ d=k$.
 $$
 \xymatrix {
 U
 \ar@/_3ex/@{-->}[rrrr]_-{w} &&
 T
 \ar[ll]_-{d}
 \ar[rr]^-{k} &&
 R
 }
 $$
 \end {claim}

 \begin {demo} [Proof.]
 Since $d_m\:T_m\to U_m$ is injective,
 there exists an archism $v\:U_m\to R_m$ such that
 $v\circ d_m=k_m$.
 By Corollary~3.2,
 $v$ is $r$-gentle for some $r\in\N$.
 Take $q\in\N$.
 We have the commutative diagram
 $$
 \xymatrix {
 U_q
 \ar[d]_-{U(m,q)} &&
 T_q
 \ar[ll]_-{d_q}
 \ar[rr]^-{k_q}
 \ar[d]^-{T(m,q)} &&
 R_q
 \ar[d]^-{R(m,q)} \\
 U_m^{[m|q]}
 \ar@/_3ex/[rrrr]_-{v^{[m|q]}} &&
 T_m^{[m|q]}
 \ar[ll]_-{d_m^{[m|q]}}
 \ar[rr]^-{k_m^{[m|q]}} &&
 R_m^{[m|q]}.
 }
 $$
 By Lemma~3.5,
 the archism $v^{[m|q]}$ is $r$-gentle.
 Since the $\Z_p$-homomorphism $R(m,q)$ is injective,
 there exists a $\Z_p$-homomorphism $f\:R_m^{[m|q]}\to R_q$
 such that
 $f\circ R(m,q)=\id$.
 Consider the $r$-gentle archism
 $$
 \xymatrix {
 w_q\:
 U_q
 \ar[rr]^-{U(m,q)} &&
 U_m^{[m|q]}
 \ar[rr]^-{v^{[m|q]}} &&
 R_m^{[m|q]}
 \ar[rr]^-{f} &&
 R_q.
 }
 $$
 Using the diagram,
 we get $w_q\circ d_q=k_q$.
 \qed
 \end {demo}


 \begin {claim} [5.2. Lemma.]
 Let
 $M$ be a simplicial crew,
 $U$ and $V$ be simplicial abeliam groups, and
 $t\:U\qto V$ be an $r$-gentle ($r\in\N$) quasisimplicial
 archism.
 Then the archism $t_\$\:\qAr(M,U)\to\qAr(M,V)$,
 $f\mapsto t\circ f$, is $r$-gentle.
 \end {claim}

 \begin {demo} [Proof.]
 This follows from Lemma~3.5
 because of the commutative diagram
 $$
 \xymatrix {
 \qAr(M,U)
 \ar[rrr]^-{t_\$}
 \ar@{=}[d] &&&
 \qAr(M,V)
 \ar@{=}[d] \\
 \prod\limits_{q\in\N,k\in M_q^\times}U_q
 \ar[rrr]^-{\prod\limits_{q\in\N,k\in M_q^\times}t_q} &&&
 \prod\limits_{q\in\N,k\in M_q^\times}V_q,
 }
 $$
 where
 $M_q^{\times}=M_q\setminus\{\textrm{basepoint}\}$.
 \qed
 \end {demo}

 \begin {claim} [5.3. Lemma.]
 Let
 $M$ and $T$ be simplicial crews,
 $U$ and $R$ be simplicial $\Z_p$-modules,
 $d\:T\to U$ and $k\:T\to R$ be simplicial archisms, and
 $w\:U\qto R$ be an $r$-gentle ($r\in\N$) quasisimplicial
 archism such that
 $w\circ d=k$.
 Then there exists an $r$-gentle quasisimplicial archism
 $z\:U^M\qto R^M$ such that
 $z\circ d^M=k^M$.
 $$
 \xymatrix {
 U^M
 \ar@/_3ex/@{-->}[rrrr]^-{z} &&
 T^M
 \ar[ll]_-{d^M}
 \ar[rr]^-{k^M} &&
 R^M
 }
 $$
 \end {claim}

 \begin {demo} [Proof.]
 Take $q\in\N$.
 We have the commutative diagram
 $$
 \xymatrix {
 (U^M)_q
 \ar[d]_-{i} &&
 (T^M)_q
 \ar[ll]_-{(d^M)_q}
 \ar[rr]^-{(k^M)_q} &&
 (R^M)_q
 \ar[d]^-{j} \\
 \qAr(\Delta^q_+\smash M,U)
 \ar[rrrr]^-{w_\$} && &&
 \qAr(\Delta^q_+\smash M,R),
 }
 $$
 where
 the $\Z_p$-homomorphism
 $i\:(U^M)_q=\sAr(\Delta^q_+\smash M,U)\to
 \qAr(\Delta^q_+\smash M,U)$ is the inclusion
 and
 $j$ is analogous.
 By Lemma~5.2,
 the archism $w_\$$ is $r$-gentle.
 There is a $\Z_p$-homomorphism
 $f\:\qAr(\Delta^q_+\smash M,R)\to(R^M)_q$ such that
 $f\circ j=\id$.
 Consider the $r$-gentle archism
 $$
 \xymatrix {
 z_q\:
 (U^M)_q
 \ar[r]^-{i} &
 \qAr(\Delta^q_+\smash M,U)
 \ar[r]^-{w_\$} &
 \qAr(\Delta^q_+\smash M,R)
 \ar[r]^-{f} &
 (R^M)_q.
 }
 $$
 Using the diagram,
 we get $z_q\circ(d^M)_q=(k^M)_q$.
 \qed
 \end {demo}


 \head {\S~6. Harmonic cofibrations}


 Let
 $T$ be a simplicial crew and
 $U$ be a simplicial abelian group.
 A cofibration $d\:T\to U$ is called {\it $r$-harmonic\/}
 ($r\in\N$) if,
 for
 any compact simplicial crews $L$ and $M$ and
 any isotypical cofibration $j\:L\to M$,
 there exist
 a simplicial archism $x\:T^L\to T^M$ and
 an $r$-gentle quasisimplicial archism $y\:U^L\qto U^M$
 such that
 $d^M\circ x=y\circ d^L$ and
 $T^j\circ x=\id$.
 $$
 \xymatrix {
 T^L
 \ar@/^2ex/[rr]^-{x}
 \ar[d]_-{d^L} &&
 T^M
 \ar[ll]^-{T^j}
 \ar[d]^-{d^M} \\
 U^L
 \ar@/_2ex/@{-->}[rr]_-{y} &&
 U^M
 \ar[ll]_-{U^j}
 }
 $$
 A cofibration is {\it harmonic\/} if
 it is $r$-harmonic for some $r\in\N$.

 By the {\it height\/} of a 0-connected space $Y$ we mean the
 supremum of those $q\in\N$ for which $\pi_q(Y)\ne1$
 (the supremum of the empty set is 0).

 \begin {claim} [6.1. Lemma.]
 Let
 $p$ be a prime and
 $Y$ be a connected CW-complex
 of finite height
 with $p$-finite homotopy groups.
 Then there exist
 a gradual simplicial crew $T$ with $\[T\]\simeq Y$,
 a gradual simplicial $\Z_p$-module $U$, and
 a harmonic cofibration $d\:T\to U$.
 \end {claim}

 \begin {demo} [Proof.]
 (Induction along the Postnikov decomposition of $Y$ with
 fibres of the form $\K(\Z_p,q)$.)
 Let $n$ be the height of $Y$.
 If $n=0$,
 then
 $Y$ is contractible,
 we put $T=U=0$ and
 that is all.
 Otherwise,
 choose an order $p$ element $e\in\pi_n(Y)$ fixed by the
 canonical action of $\pi_1(Y)$.
 Its existence follows from the well-known congruence
 $\#{\Fix_GX}\equiv\#X\pmod p$ for an action of a $p$-finite
 group $G$ on a finite set $X$
 (cf.\ the remark in \cite[Ch.~II, Example~5.2(iv)]{BouKan}).
 We attach cells to $Y$ to get a map $Y\to\=Y$
 inducing
 isomorphisms on $\pi_q$, $q\neq n$, and
 an epimorphism with the kernel generated by $e$ on $\pi_n$.
 The space $Y$ is homotopy equivalent to the homotopy fibre of
 some map $\=Y\to\K(\Z_p,n+1)$ \cite[Lemma~4.70]{Hat}.

 We assume (as an induction hypothesis) that
 there are
 gradual simplicial crew $\=T$ with $\[\=T\]\simeq\=Y$,
 gradual simplicial $\Z_p$-module $\=U$, and
 $r$-harmonic ($r\ge1$) cofibration $\=d\:\=T\to\=U$.

 Let
 $R$ be a gradual $(n+1)$-soluble simplicial $\Z_p$-module
 with $\[R\]\simeq\K(\Z_p,n+1)$,
 $Q$ be a gradual simplicial $\Z_p$-module isotypic to a point,
 and
 $c\:Q\to R$ be a fibring simplicial $\Z_p$-homomorphism
 (see \cite{Dou}).
 There is a Cartesian square of simplicial crews and archisms
 $$
 \xymatrix {
 T
 \ar[rr]^-{h}
 \ar[d]_-{f} &&
 Q
 \ar[d]^-{c} \\
 \=T
 \ar[rr]^-{k} &&
 R,
 }
 $$
 where $\[T\]\simeq Y$.
 Put $U=\=U\times Q$.
 Let $\Z_p$-homomorphisms $a\:U\to\=U$ and $b\:U\to Q$ be the
 projections.
 Let $d\:T\to U$ be the simplicial archism given by the
 conditions
 $a\circ d=\=d\circ f$ and
 $b\circ d=h$.
 Obviously, $d$ is a cofibration.

 By Lemma~5.1,
 for some $s\ge1$
 there is an $s$-gentle quasisimplicial archism $w\:\=U\qto R$
 such that
 $w\circ\=d=k$.

 We show that
 $d$ is $rs$-harmonic.
 Take
 compact simplicial crews $L$ and $M$ and
 an isotypical cofibration $j\:L\to M$.
 We need
 a simplicial archism $x\:T^L\to T^M$ and
 an $rs$-gentle quasisimplicial archism $y\:U^L\qto U^M$
 such that
 $d^M\circ x=y\circ d^L$ and
 $T^j\circ x=\id$.
 Since $\=d$ is $r$-harmonic,
 there are
 a simplicial archism $\=x\:\=T^L\to\=T^M$ and
 an $r$-gentle quasisimplicial archism $\=y\:\=U^L\qto\=U^M$
 such that
 $\=d^M\circ\=x=\=y\circ\=d^L$ and
 $\=T^j\circ\=x=\id$.

 We have the commutative square of simplicial $\Z_p$-modules
 and $\Z_p$-homomorphisms with a half-key
 $$
 \xymatrix {
 Q^L
 \ar@/^2ex/[rr]^-{t'}
 \ar[d]_-{c^L} &&
 Q^M
 \ar[ll]^-{Q^j}
 \ar[d]_-{c^M} \\
 R^L &&
 R^M
 \ar[ll]_-{R^j}
 \ar@/_2ex/[u]_-{t''}
 }
 $$
 (the half-key exists by Lemma~4.6).
 We have the simplicial archism
 $$
 u=t'\circ h^L+t''\circ k^M\circ\=x\circ f^L\:T^L\to Q^M.
 $$
 We have $c^L\circ h^L=k^L\circ f^L=k^L\circ\=
 T^j\circ\=x\circ f^L=R^j\circ k^M\circ\=x\circ f^L$.
 Therefore,
 by Lemma~4.1,
 $Q^j\circ u=h^L$ and
 $c^M\circ u=k^M\circ\=x\circ f^L$.

 Define the desired $x$ by the conditions
 $f^M\circ x=\=x\circ f^L$ and
 $h^M\circ x=u$:
 $$
 \xymatrix {
 T^M
 \ar[rr]^-{h^M}
 \ar[dd]_-{f^M} & &
 Q^M
 \ar[dd]^-{c^M} \\
 &
 T^L
 \ar[ul]_-{x}
 \ar[ur]^-{u}
 \ar[dl]_-{\=x\circ f^L} &
 \\
 \=T^M
 \ar[rr]^-{k^M} & &
 R^M.
 }
 $$
 This is possible because
 the square is Cartesian and
 the conditions are compatible:
 $k^M\circ\=x\circ f^L=c^M\circ u$.
 We have $T^j\circ x=\id$ because
 $f^L\circ T^j\circ x
 =\=T^j\circ f^M\circ x
 =\=T^j\circ\=x\circ f^L
 =f^L$ and
 $h^L\circ T^j\circ x
 =Q^j\circ h^M\circ x
 =Q^j\circ u
 =h^L$.

 By Lemma~5.3,
 there is an $s$-gentle quasisimplicial archism
 $z\:\=U^M\qto R^M$ such that
 $z\circ\=d^M=k^M$.
 We have the quasisimplicial archism
 $$
 v=t'\circ b^L+t''\circ z\circ\=y\circ a^L\:U^L\qto Q^M.
 $$
 By Lemma~3.3,
 it is $rs$-gentle.

 Define the desired $y$ by the conditions
 $a^M\circ y=\=y\circ a^L$ and
 $b^M\circ y=v$:
 $$
 \xymatrix {
 U^M
 \ar[rr]^-{b^M}
 \ar[dd]_-{a^M} & &
 Q^M \\
 &
 U^L
 \ar@{-->}[ul]_-{y}
 \ar@{-->}[ur]^-{v}
 \ar@{-->}[dl]_-{\=y\circ a^L} &
 \\
 \=U^M. & &
 }
 $$
 This is possible because
 $(a^M,b^M)\:U^M\to\=U^M\times Q^M$ is an isomorphism.
 Obviously,
 $y$ is $rs$-gentle.
 We have $d^M\circ x=y\circ d^L$ because
 $a^M\circ d^M\circ x
 =\=d^M\circ f^M\circ x
 =\=d^M\circ\=x\circ f^L
 =\=y\circ\=d^L\circ f^L
 =\=y\circ a^L\circ d^L
 =a^M\circ y\circ d^L$ and
 $b^M\circ d^M\circ x
 =h^M\circ x
 =u
 =t'\circ h^L+t''\circ k^M\circ\=x\circ f^L
 =t'\circ h^L+t''\circ z\circ\=d^M\circ\=x\circ f^L
 =t'\circ h^L+t''\circ z\circ\=y\circ\=d^L\circ f^L
 =t'\circ b^L\circ d^L+t''\circ z\circ\=y\circ a^L\circ d^L
 =v\circ d^L
 =b^M\circ y\circ d^L$.

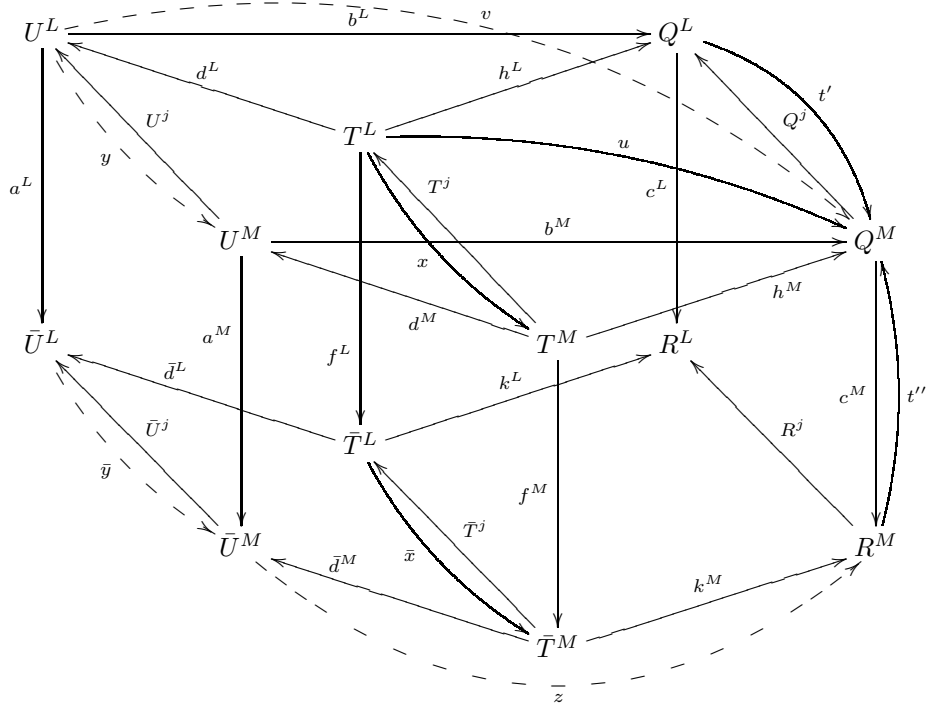
 $$
 \xymatrix {
 U^L
 \ar[ddd]_-{a^L}
 \ar@/_2ex/@{-->}[ddrr]_-{y}
 \ar[rrrrrr]^-{b^L}
 \ar@/^10ex/@{-->}[ddrrrrrrrr]^-{v} & & & & & &
 Q^L
 \ar[ddd]_-{c^L}
 \ar@/^4ex/[ddrr]^-{t'} & &
 \\
 & & &
 T^L
 \ar[ulll]_-{d^L}
 \ar[ddd]_(0.7){f^L}
 \ar@/_2ex/[ddrr]_-{x}
 \ar[urrr]^-{h^L}
 \ar@/^3ex/[drrrrr]^-{u} & & & & &
 \\
 & &
 U^M
 \ar[uull]_-{U^j}
 \ar[ddd]_(0.3){a^M}
 \ar[rrrrrr]^-{b^M} & & & & & &
 Q^M
 \ar[uull]_-{Q^j}
 \ar[ddd]_-{c^M}
 \\
 \=U^L
 \ar@/_2ex/@{-->}[ddrr]_-{\=y} & & & & &
 T^M
 \ar[uull]_(0.7){T^j}
 \ar[ulll]^(0.4){d^M}
 \ar[ddd]_-{f^M}
 \ar[urrr]_(0.7){h^M} &
 R^L & &
 \\
 & & &
 \=T^L
 \ar[ulll]_(0.6){\=d^L}
 \ar@/_2ex/[ddrr]_(0.4){\=x}
 \ar[urrr]^-{k^L} & & & & &
 \\
 & &
 \=U^M
 \ar[uull]_-{\=U^j}
 \ar@/_12ex/@{-->}[rrrrrr]_-{z} & & & & & &
 R^M
 \ar@/_2ex/[uuu]_-{t''}
 \ar[uull]_-{R^j}
 \\
 & & & & &
 \=T^M
 \ar[uull]_-{\=T^j}
 \ar[ulll]_(0.7){\=d^M}
 \ar[urrr]^-{k^M} & & &
 }
 $$
 (The straight arrows of this diagram form a commutative
 subdiagram.)
 \qed
 \end {demo}


 \head {\S~7. Two filtrations of the module $C_0(U^K)$}


 \begin {claim} [7.1. Lemma.]
 Let $U_i$, $i\in I$, be a finite collection of abelian groups.
 Put
 $$
 U_J=\bigoplus_{i\in J}U_i,
 \qquad
 J\subseteq I,
 $$
 and
 $U=U_I$.
 Let $p_J\:U\to U_J$ be the projections.
 Then
 for any $r\in\N$
 $$
 \bigcap_{J\subseteq I\,:\,\#J\le r}
 \ker(p_J)_\R
 \subseteq\ai\R U^{r+1}.
 $$
 in the $\R$-algebra $\ga\R U$.
 \end {claim}

 \begin {demo} [Proof.]
 Let $s_J\:U_J\to U$ be the canonical embeddings.
 Put $q_J=s_J\circ p_J\:U\to U$.
 We assume $\#I>r$
 (otherwise, the assertion is trivial).
 For $u\in U$,
 we have
 (cf.\ \cite[Lemma~5.5]{Gou})
 \begin {multline*}
 \`u\'-
 \sum_{J\subseteq I\,:\,\#J\le r}
 (-1)^{r-\#J}\bc{\#I-\#J-1}{r-\#J}\`q_J(u)\'= \\ =
 \sum_{J\subseteq I}
 \big(
 \sum_{M\subseteq I\,:\,M\supseteq J,\,\#M>r}(-1)^{\#M-\#J}
 \big)
 \`q_J(u)\'= \\ =
 \sum_{M\subseteq I\,:\,\#M>r}
 \big(\sum_{J\subseteq M}(-1)^{\#M-\#J}\`q_J(u)\'
 \big)= \\ =
 \sum_{M\subseteq I\,:\,\#M>r}
 \prod_{i\in M}
 (\`q_{\{i\}}(u)\'-1)\in
 \ai\R U^{r+1}.
 \end {multline*}
 It follows that
 for $w\in\ga\R U$
 we have
 $$
 w-
 \sum_{J\subseteq I\,:\,\#J\le r}
 (-1)^{r-\#J}\bc{\#I-\#J-1}{r-\#J}(q_J)_\R(w)\in
 \ai\R U^{r+1}.
 $$
 If
 $$
 w\in
 \bigcap_{J\subseteq I\,:\,\#J\le r}
 \ker(p_J)_\R,
 $$
 then,
 using that $\ker(p_J)_\R=\ker(q_J)_\R$,
 we get $w\in\ai\R U^{r+1}$.
 \qed
 \end {demo}


 For a simplicial abelian group $V$,
 the module $C_0(V)=\ga\R{V_0}$
 has the filtration $C_0^{\%s}(V)=\ai\R{V_0}^s$, $s\in\N$.

 \begin {claim} [7.2. Corollary.]
 Let
 $K$ be a compact simplicial crew,
 $E$ be a field,
 $U$ be a simplicial $E$-module, and
 $r\in\N$ be a number.
 Consider the $\R$-homomorphism
 $$
 \xymatrix {
 C_0(U^K)
 \ar[rr]^-{\|KU\mu_r} &&
 \Hom_0(C_*(K^r),C_*(U^r)).
 }
 $$
 Then $\ker\|KU\mu_r\subseteq C_0^{\%r+1}(U^K)$.
 \end {claim}

 \begin {demo} [Proof.]
 Take an element $B\in\ker\|KU\mu_r$.
 We show that $B\in C_0^{\%r+1}(U^K)$.

 There is $n\in\N$ such that
 the simplicial crew $K$ is generated by a finite collection of
 $n$-simplices: $g_i\in K_n$, $i\in I$.
 We have the $E$-homomorphism $h\:(U^K)_0\to U_n^I$,
 $b\mapsto(b(g_i))_{i\in I}$.
 It is injective.
 Therefore,
 there is
 an $E$-homomorphism $f\:U_n^I\to(U^K)_0$ such that
 $f\circ h=\id$.
 It suffices to show that
 $h_\R(B)\in\ai\R{U_n^I}^{r+1}$.
 Indeed, then
 $B
 =f_\R(h_\R(B))
 \in\ai\R{(U^K)_0}^{r+1}
 =C_0^{\%r+1}(U^K)$.

 For $J\subseteq I$,
 let $p_J\:U_n^I\to U_n^J$ be the projection.
 Take $J\subseteq I$ with $\#J\le r$.
 By Lemma~7.1,
 it suffices to verify that
 $(p_J)_\R(h_\R(B))=0$.

 Choose
 a function $t\:J\to\{1,\dotsc,r\}$ and
 a simplex $k=(k_1,\dotsc,k_r)\in K_n^r$
 such that
 $k_{t(i)}=g_i$, $i\in J$.
 We have
 the $E$-homomorphism $U_n^t\:U_n^r\to U_n^J$,
 the $\R$-homomorphism
 $(U_n^t)_\R\:C_n(U^r)=\ga\R{U_n^r}\to\ga\R{U_n^J}$, and
 the commutative diagram
 $$
 \xymatrix {
 \ga\R{(U^K)_0}
 \ar[rrr]^-{h_\R}
 \ar[d]_-{\|KU\mu_r} &&&
 \ga\R{U_n^I}
 \ar[d]^-{(p_J)_\R} \\
 \Hom_0(C_*(K^r),C_*(U^r))
 \ar[rrr]^-{v\mapsto(U_n^t)_\R(v(\`k\'))} &&&
 \ga\R{U_n^J}.
 }
 $$
 Since $\|KU\mu_r(B)=0$,
 we get $(p_J)_\R(h_\R(B))=0$.
 \qed
 \end {demo}


 \head {\S~8. Simplicial approximation}


 \begin {claim} [8.1. Lemma.]
 Let
 $K$ be a compact simplicial crew,
 $W$ be a simplicial crew, and
 $f\:\[K\]\to\[W\]$ be a map.
 Then there exist
 a compact simplicial crew $L$,
 an isotypy $e\:L\to K$, and
 a simplicial archism $g\:L\to W$
 such that
 $f\circ\[e\]\sim\[g\]$.
 \end {claim}

 \begin {demo}
 See \cite[Corollary~4.8]{Jar}.
 \qed
 \end {demo}


 For simplicial crews $L$ and $T$,
 the geometrical realization $\[?\]\:(T^L)_0\to\[T\]^{\[L\]}$
 induces an $\R$-homomorphism
 $\(?\)\:H_0(T^L)\to H_0(\[T\]^{\[L\]})$.

 \begin {claim} [8.2. Lemma.]
 Let
 $K$ be a compact simplicial crew,
 $T$ be a simplicial crew, and
 $r\in\N$ be a number.
 Then,
 for any $A\in\ker\|{\[K\]}{\[T\]}\mu_r$,
 there exist
 a compact simplicial crew $L$,
 an isotypy $e\:L\to K$, and
 an element $B\in\ker\|LT\mu_r$ such that
 $H_0(\[T\]^{\[e\]})(\|{\[K\]}{\[T\]}\nu(A))=\(\|LT\nu(B)\)$:
 $$
 \xymatrix {
 \Hom_0(C_*(L^r),C_*(T^r)) &&
 C_0(T^L)
 \ar@{:=}@(dl,dr)[]^-{B}
 \ar[ll]_-{\|LT\mu_r}
 \ar[rr]^-{\|LT\nu} &&
 H_0(T^L)
 \ar[d]^-{\(?\)} \\
 && &&
 H_0(\[T\]^{\[L\]}) \\
 \Hom(C_0(\[K\]^r),C_0(\[T\]^r)) &&
 C_0(\[T\]^{\[K\]})
 \ar@{:=}@(ul,ur)[]_-{A}
 \ar[ll]_-{\|{\[K\]}{\[T\]}\mu_r}
 \ar[rr]^-{\|{\[K\]}{\[T\]}\nu} &&
 H_0(\[T\]^{\[K\]}).
 \ar[u]_-{H_0(\[T\]^{\[e\]})}
 }
 $$
 \end {claim}

 \begin {demo} [Proof.]
 We have
 $$
 A=\sum_{i=1}^mv_i\`a_i\',
 $$
 where
 $m\in\N$,
 $v_i\in\R$, and
 $a_i\in\[T\]^{\[K\]}$.
 For $x\in\[K\]$,
 define an equivalence (relation) $c(x)$ on the set
 $I=\{1,\dotsc,m\}$:
 $c(x)=\{\,(i,j):a_i(x)=a_j(x)\,\}$.
 Put $E=\{\,c(x):x\in\[K\]\,\}$.

 We call an equivalence on $I$ {\it neutral\/} if
 $$
 \sum_{i\in J}v_i=0
 $$
 for all its classes $J\subseteq I$.
 We show that
 for any $h_1,\dotsc,h_r\in E$
 the equivalence $h=h_1\cap\dotso\cap h_r$ is neutral.
 For each $s=1,\dotsc,r$,
 there is a point $x_s\in\[K\]$ such that
 $h_s=c(x_s)$.
 Put $x=(x_1,\dotsc,x_r)\in\[K\]^r$.
 In $C_0(\[T\]^r)$,
 we have
 $$
 \sum_{i\in I}
 v_i\`a_i^r(x)\'=
 \|{\[K\]}{\[T\]}\mu_r(A)=0.
 $$
 It follows that
 $h$ is neutral
 because
 $$
 a_i^r(x)=a_j^r(x)
 \iff
 (i,j)\in h
 $$
 for $i,j\in I$.

 For each equivalence $h$ on $I$,
 there is the corresponding simplicial subcrew
 $V(h)\subseteq T^m$ (the diagonal):
 $$
 V(h)_q=
 \{\,(t_1,\dotsc,t_m)\in T^m_q:
 \textrm{$t_i=t_j$ for all $(i,j)\in h$}\,\}.
 $$
 Put
 $$
 W=\bigcup_{h\in E}V(h)\subseteq T^m.
 $$
 We have the maps
 $a=(a_1,\dotsc,a_m)\:\[K\]\to\[T\]^m$ and
 $\~a=d^{-1}\circ a\:\[K\]\to\[T^m\]$,
 where $d\:\[T^m\]\to\[T\]^m$ is the canonical bijective map.
 For $x\in\[K\]$,
 we have $\~a(x)\in\[V(c(x))\]$.
 Therefore $\im\~a\subseteq\[W\]$.
 Using Lemma~8.1,
 we find
 a compact simplicial crew $L$,
 an isotypy $e\:L\to K$, and
 a simplicial archism $b=(b_1,\dotsc,b_m)\:L\to T^m$
 such that
 $\im b\subseteq W$ and
 $\~a\circ\[e\]\sim\[b\]$.
 Put
 $$
 B=\sum_{i=1}^mv_i\`b_i\'.
 $$
 We have $a_i\circ\[e\]\sim\[b_i\]$.
 Therefore $H_0(\[T\]^{\[e\]})(\|{\[K\]}{\[T\]}\nu(A))=
 \(\|LT\nu(B)\)$.
 We show that
 $\|KT\mu_r(B)=0$.
 For $k=(k_1,\dotsc, k_r)\in K^r_q$ ($q\in\N$),
 we have
 $$
 \|KT\mu_r(B)(\`k\')=
 \sum_{i=1}^mv_i\`b_i^r(k)\'.
 $$
 Take $s=1,\dotsc,r$.
 Since $\im b\subseteq W$,
 there is $h_s\in E$ such that
 $b(k_s)\in V(h_s)$.
 Therefore,
 the function $i\mapsto b_i(k_s)$ is subordinate to
 (i.~e.\ constant on the classes of)
 the equivalence $h_s$.
 Since $b_i^r(k)=(b_i(k_1),\dotsc,b_i(k_r))$,
 the function $i\mapsto b_i^r(k)$ is subordinate to the
 equivalence $h=h_1\cap\dotso\cap h_r$.
 Since $h$ is neutral,
 we get $\|KT\mu_r(B)(\`k\')=0$.
 \qed
 \end {demo}


 \head {\S~9. The inclusion $\ker\|XY\mu_r\subseteq\ker\|XY\nu$
        for large $r$}


 \begin {claim} [9.1. Lemma.]
 Let $X$, $Y$, $\~X$, and $\~Y$ be spaces.
 Suppose that $X\simeq\~X$ and $Y\simeq\~Y$.
 Then,
 for any $r\in\N$,
 we have
 $$
 \ker\|XY\mu_r\subseteq\ker\|XY\nu
 \iff
 \ker\|{\~X}{\~Y}\mu_r\subseteq\ker\|{\~X}{\~Y}\nu.
 $$
 \end {claim}

 \begin {demo} [Proof.]
 There are homotopy euivalences
 $k\:X\to\~X$ and
 $h\:\~Y\to Y$.
 We have the commutative diagram of $\R$-modules and
 $\R$-homomorphisms:
 $$
 \xymatrix {
 \Hom(C_0(\~X^r),C_0(\~Y^r))
 \ar[d] &&
 C_0(\~Y^{\~X})
 \ar[ll]_-{\|{\~X}{\~Y}\mu_r}
 \ar[rr]^-{\|{\~X}{\~Y}\nu}
 \ar[d]^-{C_0(h^k)} &&
 H_0(\~Y^{\~X})
 \ar[d]^-{H_0(h^k)} \\
 \Hom(C_0(X^r),C_0(Y^r)) &&
 C_0(Y^X)
 \ar[ll]_-{\|XY\mu_r}
 \ar[rr]^-{\|XY\nu} &&
 H_0(Y^X),
 }
 $$
 where the vertical arrows are induced by $k$ and $h$.
 Since $H_0(h^k)$ is an isomorphism,
 we get the implication $\Rightarrow$.
 The implication $\Leftarrow$ is analogous.
 \qed
 \end {demo}


 Let $p$ be a prime.
 Assume $\R=\Z_p$.

 \begin {claim} [9.2.]
 Let
 $X$ be a compact CW-complex and
 $Y$ be a connected CW-complex
 of finite height
 with $p$-finite homotopy groups.
 Then,
 for any sufficiently large $r\in\N$,
 we have $\ker\|XY\mu_r\subseteq\ker\|XY\nu$
 in the diagram
 $$
 \xymatrix {
 \Hom(C_0(X^r),C_0(Y^r)) &&
 C_0(Y^X)
 \ar[ll]_-{\|XY\mu_r}
 \ar[rr]^-{\|XY\nu} &&
 H_0(Y^X).
 }
 $$
 \end {claim}

 \begin {demo} [Proof.]
 By Lemma~6.1,
 for some $s\in\N$,
 there are
 a gradual simplicial crew $T$ with $\[T\]\simeq Y$,
 a gradual simplicial $\Z_p$-module $U$, and
 an $s$-harmonic cofibration $d\:T\to U$.
 We have $X\simeq\[K\]$
 for some compact simplicial crew $K$.
 Obviously,
 $(U^K)_0$ is a finite $\Z_p$-module.
 By Lemma~3.1,
 $C_0^{\%t+1}(U^K)=0$
 for some $t\in\N$.
 Take $r\ge st$.
 We show that
 $\ker\|{\[K\]}{\[T\]}\mu_r\subseteq\ker\|{\[K\]}{\[T\]}\nu$
 in the diagram
 $$
 \xymatrix {
 \Hom(C_0(\[K\]^r),C_0(\[T\]^r)) &&
 C_0(\[T\]^{\[K\]})
 \ar[ll]_-{\|{\[K\]}{\[T\]}\mu_r}
 \ar[rr]^-{\|{\[K\]}{\[T\]}\nu} &&
 H_0(\[T\]^{\[K\]}).
 }
 $$
 This will suffice
 by Lemma~9.1.

 Take an element $A\in\ker\|{\[K\]}{\[T\]}\mu_r$.
 We show that $A\in\ker\|{\[K\]}{\[T\]}\nu$.
 By Lemma~8.2,
 there are
 a compact simplicial crew $L$,
 an isotypy $e\:L\to K$, and
 an element $B\in\ker\|LT\mu_r$ such that
 $H_0(\[T\]^{\[e\]})(\|{\[K\]}{\[T\]}\nu(A))=\(\|LT\nu(B)\)$.
 Since $\[e\]$ is a homotopy equivalence,
 $H_0(\[T\]^{\[e\]})$ is an isomorphism.
 Therefore
 it suffices to show that
 $\|LT\nu(B)=0$.

 Let a simplicial crew $M$ be the (reduced) cylinder of $e$.
 We have the homotopy commutative diagram
 $$
 \xymatrix {
 K
 \ar[dr]_-{i} & &
 L
 \ar[ll]_-{e}
 \ar[dl]^-{j} \\
 &
 M, &
 }
 $$
 where $i$ and $j$ are the canonical cofibrations.
 By the definition of a cylinder,
 $i$ is an isotypy.
 Since $e$ is an isotypy,
 $j$ is an isotypy too.
 Since $d$ is $s$-harmonic,
 there is the commutative diagram
 $$
 \xymatrix {
 T^L
 \ar@/^3ex/[rrrr]^-{\id}
 \ar[rr]_-{x}
 \ar[d]_-{d^L} &&
 T^M
 \ar[rr]_-{T^j}
 \ar[d]^-{d^M} &&
 T^L \\
 U^L
 \ar@{-->}[rr]^-{y} &&
 U^M, &&
 }
 $$
 where
 $x$ is a simplicial archism and
 $y$ is an $s$-gentle quasisimplicial archism.
 We have the commutative diagram of $\Z_p$-homomorphisms:
 $$
 \xymatrix {
 \Hom_0(C_*(L^r),C_*(T^r))
 \ar[d] &
 C_0(T^L)
 \ar@{:=}@(ul,ur)[]_-{B}
 \ar[l]_-{\|LT\mu_r}
 \ar[r]^-{C_0(x)}
 \ar[d]^-{C_0(d^L)} &
 C_0(T^M)
 \ar@{:=}@(ul,ur)[]_-{B_1}
 \ar[r]^{C_0(T^i)}
 \ar[d]^-{C_0(d^M)} &
 C_0(T^K)
 \ar@{:=}@(ul,ur)[]_-{B_2}
 \ar[d]^-{C_0(d^K)} \\
 \Hom_0(C_*(L^r),C_*(U^r)) &
 C_0(U^L)
 \ar@{:=}@(dl,dr)[]^-{B'}
 \ar[l]_-{\|LU\mu_r}
 \ar[r]^-{C_0(y)} &
 C_0(U^M)
 \ar@{:=}@(dl,dr)[]^-{B'_1}
 \ar[r]^{C_0(U^i)} &
 C_0(U^K)
 \ar@{:=}@(dl,dr)[]^-{B'_2},
 }
 $$
 where the vertical arrows are induced by the cofibration $d$;
 $B_1,\dotsc,B'_2$ are the images of $B$ in the corresponding
 modules.
 Since $\|LT\mu_r(B)=0$,
 we have $\|LU\mu_r(B')=0$.
 By Corollary~7.2,
 $B'\in C_0^{\%r+1}(U^L)$.
 Since
 $r\ge st$
 and the archism $y_0$ is $s$-gentle,
 we have,
 by Corollary~3.4,
 $B'_1\in C_0^{\%t+1}(U^M)$.
 Since $(U^i)_0$ is a homomorphism,
 $B'_2\in C_0^{\%t+1}(U^K)$.
 We have $C_0^{\%t+1}(U^K)=0$.
 It follows that
 $B'_2=0$.
 Since $d$ is a cofibration,
 $C_0(d^K)$ is injective.
 Therefore
 $B_2=0$.

 We have the commutative diagram of $\Z_p$-homomorphisms
 $$
 \xymatrix {
 C_0(T^L)
 \ar@{:=}@(dl,dr)[]^-{B}
 \ar[rr]^-{\id}
 \ar[dr]_-{C_0(x)} & &
 C_0(T^L)
 \ar@{:=}@(ul,ur)_-{B}
 \ar[rr]^-{\|LT\nu} & &
 H_0(T^L) \\
 &
 C_0(T^M)
 \ar@{:=}@(dl,dr)[]^-{B_1}
 \ar[rr]^-{\|MT\nu}
 \ar[ur]^-{C_0(T^j)}
 \ar[dr]_-{C_0(T^i)} & &
 H_0(T^M)
 \ar[ur]^-{H_0(T^j)}
 \ar[dr]_-{H_0(T^i)} &
 \\
 & &
 C_0(T^K)
 \ar@{:=}@(dl,dr)[]^-{B_2}
 \ar[rr]^-{\|KT\nu} & &
 H_0(T^K)
 \ar[uu]_-{H_0(T^e)}.
 }
 $$
 Since $B_2=0$,
 we get $\|LT\nu(B)=0$.
 \qed
 \end {demo}


 Consider the filtration of the complex $C_*(Y^X)$ formed by
 the kernels of the $\Z_p$-homomorphisms
 $$
 \xymatrix {
 C_q(Y^X)
 \ar[r]^-{i_q} &
 C_0(Y^{\Delta^q_+\smash X})
 \ar[rr]^-{\displaystyle\|{\Delta^q_+\smash X}Y\mu_r} &&
 \Hom(C_0((\Delta^q_+\smash X)^r),C_0(Y^r)),
 }
 $$
 where $i_q$ are the obvious isomorphisms.
 Does this filtration converge?


 \head {\S~10. Deducing Theorem~1.2 from claim~9.2}


 \begin {claim} [10.1. Lemma.]
 Let
 $X$, $Y$, $\~X$, and $\~Y$ be spaces,
 $k\:X\to\~X$ and $h\:\~Y\to Y$ be maps,
 $V$ be an abelian group, and
 $f\:[X,Y]\to V$ be an invariant.
 Consider the invariant $\~f\:[\~X,\~Y]\to V$,
 $\~u\mapsto f([h]\circ\~u\circ[k])$.
 Then $\Deg\~f\le\Deg f$.
 \end {claim}

 \begin {demo} [Proof.]
 Take $r\in\N$.
 The maps $k$ and $h$ induce a homomorphism
 $$
 t\:
 \Hom(C_0(\~X^r),C_0(\~Y^r))\to
 \Hom(C_0(X^r),C_0(Y^r)).
 $$
 We have $t(C_0(\~a^r))=C_0((h\circ\~a\circ k)^r)$,
 $\~a\in\~Y^{\~X}$.
 Assume that $\Deg f\le r$.
 There is a homomorphism $l\:\Hom(C_0(X^r),C_0(Y^r))\to V$
 such that
 $f([a])=l(C_0(a^r))$ for all $a\in Y^X$.
 Consider the homomorphism
 $\~l=l\circ t\:\Hom(C_0(\~X^r),C_0(\~Y^r))\to V$.
 For $\~a\in\~Y^{\~X}$
 we have
 $\~f([\~a])
 =f([h\circ\~a\circ k])
 =l(C_0((h\circ\~a\circ k)^r))
 =l(t(C_0(\~a^r)))
 =\~l(C_0(\~a^r))$.
 Therefore
 $\Deg\~f\le r$.
 \qed
 \end {demo}


 \subhead {Proof of Theorem~1.2.}
 (1)
 {\it Case of $Y$ of finite height.\/}
 It suffices to show that
 the ``universal'' invariant $F\:[X,Y]\to H_0(Y^X;\Z_p)$,
 $u\mapsto\`u\'$, has finite degree.
 For $r\in\N$
 we have the commutative diagram
 $$
 \xymatrix {
 \Hom_\Z(C_0(X^r;\Z),C_0(Y^r;\Z))
 \ar[d]_-{m'} &
 C_0(Y^X;\Z)
 \ar[l]_-{\|XY{\~\mu}_r}
 \ar[d]^-{m} &
 \\
 \Hom_{\Z_p}(C_0(X^r;\Z_p),C_0(Y^r;\Z_p)) &
 C_0(Y^X;\Z_p)
 \ar[l]_-{\|XY\mu_r}
 \ar[r]^-{\|XY\nu} &
 H_0(Y^X;\Z_p),
 }
 $$
 where $m$ and $m'$ are the homomorphisms of reduction modulo
 $p$;
 the tilde over $\mu$ in the upper row means ``over $\Z$''.
 By claim~9.2,
 we have $\ker\|XY\mu_r\subseteq\ker\|XY\nu$
 for sufficiently large $r$.
 Then
 there is a $\Z_p$-homomorphism
 $t\:\Hom_{\Z_p}(C_0(X^r;\Z_p),C_0(Y^r;\Z_p))\to
 H_0(Y^X;\Z_p)$ such that
 $t\circ\|XY\mu_r=\|XY\nu$.
 For $a\in Y^X$,
 we have
 $F([a])
 =(\|XY\nu\circ m)(\`a\')
 =(t\circ m'\circ\|XY{\~\mu}_r)(\`a\')
 =(t\circ m')(C_0(a^r;\Z))$.
 Therefore
 $\Deg F\le r$.

 (2)
 {\it General case.\/}
 There are
 a connected CW-complex $\=Y$
 of finite height
 with $p$-finite homotopy groups and
 a $(\dim X+1)$-connected map $h\:Y\to\=Y$
 ($\=Y$ is obtained from $Y$ by attaching cells of high
 dimensions).
 The induced function $h_\$\:[X,Y]\to[X,\=Y]$ is bijective.
 Consider the invariant $\=f=f\circ h_\$^{-1}\:[X,\=Y]\to\Z_p$.
 By Lemma~10.1,
 $\Deg f\le\Deg\=f$.
 By (1),
 $\Deg\=f<\infty$.
 \qed


 \head {\S~11. Deducing Theorem~1.1 from Theorem~1.2}


 \begin {claim} [11.1. Lemma
                 \rm {\cite[Ch.~VI, Proposition~8.6]{BouKan}}.]
 Let
 $X$ be a connected compact CW-complex,
 $Y$ be a nilpotent connected CW-complex
 with finitely generated homotopy groups, and
 $u_1,u_2\in[X,Y]$ be distinct classes.
 Then,
 for some prime $p$,
 there exist
 a connected CW-complex $\=Y$
 with $p$-finite homotopy groups and
 a map $h\:Y\to\=Y$
 such that
 $[h]\circ u_1\ne[h]\circ u_2$ in $[X,\=Y]$.
 \qed
 \end {claim}


 \subhead {Proof of Theorem~1.1.}
 By Lemma~11.1,
 for some prime $p$
 there are
 a connected CW-complex $\=Y$
 with $p$-finite homotopy groups, and
 a map $h\:Y\to\=Y$
 such that
 the classes $\=u_i=[h]\circ u_i$, $i=1,2$, are distinct.
 There is an invariant $\=f\:[X,\=Y]\to\Z_p$ such that
 $\=f(\=u_1)\ne\=f(\=u_2)$.
 By Theorem~1.2,
 $\Deg\=f<\infty$.
 Consider the invariant $f=\=f\circ h_\$\:[X,Y]\to\Z_p$.
 By Lemma~10.1,
 $\Deg f<\infty$.
 We have
 $f(u_1)
 =\=f(\=u_1)
 \ne\=f(\=u_2)
 =f(u_2)$.
 \qed


 \head {\S~12. Properties of finite-degree invariants}


 Put $\E=\{0,1\}\subseteq\Z$.
 For $e=(e_1,\dotsc,e_n)\in\E^n$,
 put $\*e=e_1+\dotso+e_n$.

 Cosider a wedge of spaces $W=T_1\vee\dotso\vee T_n$.
 Let $\incl^W_k\:T_k\to W$ be the inclusions.
 For $e\in\E^n$,
 put $\M^W_e=m_1\vee\dotso\vee m_n\:W\to W$,
 where $m_k\:T_k\to T_k$ is:
 the identity if $e_k=1$, and
 the constant map otherwise.

 \begin {claim} [12.1. Lemma.]
 Let
 $X$ and $Y$ be spaces,
 $V$ be an abelian group,
 $f\:[X,Y]\to V$ be an invariant of degree at most $r\in\N$,
 $W=T_1\vee\dotso\vee T_{r+1}$ be a wedge of spaces, and
 $k\:X\to W$ and $h\:W\to Y$ be maps.
 Then
 $$
 \sum_{e\in\E^{r+1}}
 (-1)^{\*e}f([h\circ\M^W_e\circ k])=0.
 $$
 \end {claim}

 \begin {demo} [Proof.]
 Consider the invariant $\~f\:[W,W]\to V$,
 $\~u\mapsto f([h]\circ\~u\circ[k])$.
 We show that
 $$
 \sum_{e\in\E^{r+1}}
 (-1)^{\*e}\~f([\M^W_e])=0.
 $$
 By Lemma~10.1,
 $\Deg\~f\le r$,
 i.~e.\
 there is a homomorphism $l\:\Hom(C_0(W^r),C_0(W^r))\to V$ such
 that
 $\~f([\~a])=l(C_0(\~a^r))$ for all $\~a\in W^W$
 (hereafter, $\R=\Z$).
 Therefore
 it suffices to show that
 $$
 \sum_{e\in\E^{r+1}}
 (-1)^{\*e}C_0((\M^W_e)^r)=0.
 $$
 Take a point $w=(w_1,\dotsc,w_r)\in W^r$.
 There is $s\in\{1,\dotsc,r+1\}$ such that
 $\{w_1,\dotsc,w_r\}\cap T_s\subseteq\{\textrm{basepoint}\}$.
 The point $(\M^W_e)^r(w)\in W^r$ does not depend on the $s$th
 component of $e$.
 Since $C_0((\M^W_e)^r)(\`w\')=\`(\M^W_e)^r(w)\'$,
 it follows that
 $$
 \sum_{e\in\E^{r+1}}
 (-1)^{\*e}C_0((\M^W_e)^r)(\`w\')=0.
 $$
 \qed
 \end {demo}


 \subhead {Maps $S^n\to Y$.}
 In this subsection,
 we use multiplicative notation for homotopy groups.

 \begin {claim} [12.2. Lemma.]
 Let
 $n\ge1$ be a number,
 $Y$ be a space,
 $V$ be an abelian group, and
 $f\:\pi_n(Y)\to V$ be an invariant of degree at most $r\in\N$.
 Then $f$ is $r$-gentle.
 \end {claim}

 \begin {demo} [Proof.]
 Take elements $u_1,\dotsc,u_{r+1}\in\pi_n(Y)$.
 We show that
 $\+f((1-\`u_1\')\dotso(1-\`u_{r+1}\'))=0$.
 Put $W=S^n\vee\dotso\vee S^n$ ($r+1$ summands).
 Let
 $k\:S^n\to W$ be a map
 with $[k]=[\incl^W_1]\dotso[\incl^W_{r+1}]$ in $\pi_n(W)$, and
 $h\:W\to Y$ be a map
 with $[h\circ\incl^W_s]=u_s$ in $\pi_n(Y)$.
 By Lemma~12.1,
 $$
 \sum_{e\in\E^{r+1}}
 (-1)^{\*e}f([h\circ\M^W_e\circ k])=0.
 $$
 This is what we need because
 $[h\circ\M^W_e\circ k]=u_1^{e_1}\dotso u_{r+1}^{e_{r+1}}$ in
 $\pi_n(Y)$.
 \qed
 \end {demo}

 We denote the Whitehead product by the sign $*$.

 \begin {claim} [12.3. Lemma.]
 Let
 $m,n\ge1$ be numbers,
 $Y$ be a space, and
 $f\:\pi_{m+n-1}(Y)\to V$ be an invariant
 of degree at most $r\in\N$.
 Then the function $b\:\pi_m(Y)\times\pi_n(Y)\to V$,
 $(u,v)\mapsto f(u*v)$, is $r$-gentle.
 \end {claim}

 \begin {demo} [Proof.]
 Assume $r>0$
 (otherwise, the claim is trivial).
 Take elements
 $u_1,\dotsc,u_p\in\pi_m(Y)$ and
 $v_1,\dotsc,v_q\in\pi_n(Y)$,
 where
 $p,q\ge0$ and
 $p+q=r+1$.
 By Lemma~3.10,
 it suffices to show that
 $\+b((1-\`\^u_1\')\dotso(1-\`\^u_p\')
      (1-\`\^v_1\')\dotso(1-\`\^v_q\'))=0$,
 where
 $\^u_s=(u_s,1)\in\pi_m(Y)\times\pi_n(Y)$ and
 $\^v_s=(1,v_s)\in\pi_m(Y)\times\pi_n(Y)$.
 Put
 $W=S^m\vee\dotso\vee S^m\vee S^n\vee\dotso\vee S^n$
 ($p$ times $S^m$ and
 $q$ times $S^n$).
 Let
 $k\:S^{m+n-1}\to W$ be a map with
 $[k]=([\incl^W_1]\dotso[\incl^W_p])*
      ([\incl^W_{p+1}]\dotso[\incl^W_{r+1}])$ in
 $\pi_{m+n-1}(W)$ and
 $h\:W\to Y$ be a map with
 $[h\circ\incl^W_s]=u_s$ in $\pi_m(Y)$
 for $s=1,\dotsc,p$ and
 $[h\circ\incl^W_{p+t}]=v_t$ in $\pi_n(Y)$
 for $t=1,\dotsc,q$.
 By Lemma~12.1,
 $$
 \sum_{e\in\E^{r+1}}
 (-1)^{\*e}f([h\circ\M^W_e\circ k])=0.
 $$
 This is what we need because
 $[h\circ\M^W_e\circ k]
 =(u_1^{e_1}\dotso u_p^{e_p})
 *(v_1^{e_{p+1}}\dotso v_q^{e_{r+1}})$ in $\pi_{m+n-1}(Y)$ and,
 consequently,
 $f([h\circ\M^W_e\circ k])
 =b(u_1^{e_1}\dotso u_p^{e_p},
    v_1^{e_{p+1}}\dotso v_q^{e_{r+1}})
 =b(\^u_1^{e_1}\dotso\^u_p^{e_p}
    \^v_1^{e_{p+1}}\dotso\^v_q^{e_{r+1}})$.
 \qed
 \end {demo}


 \subhead {Maps $S^{n-1}\times S^n\to S^n_{(\Q)}$.}
 In this subsection,
 we
 prove claims 1.5--1.7 and
 use the objects defined in the corresponding subsection of
 \S~1.
 For $u\in\pi_p(Y)$ and $v\in\pi_q(Y)$,
 the class $(u,v)\in[S^p\vee S^q,Y]$ is defined in the obvious
 way.

 Let $x\:S^n\vee S^{2n-1}\to S^n\times S^{2n-1}$ be the
 canonical embedding of a wedge in the product.
 Consider the map
 $(\proj_2,c)\:S^{n-1}\times S^n\to S^n\times S^{2n-1}$,
 where
 $\proj_2\:S^{n-1}\times S^n\to S^n$ is the projection and
 $c\:S^{n-1}\times S^n\to S^{2n-1}$ is the map defined in \S~1.
 There exists a
 (unique up to homotopy)
 map $b\:S^{n-1}\times S^n\to S^n\vee S^{2n-1}$ such that
 $x\circ b\sim(\proj_2,c)$.
 For $p,q\in\Z$,
 we have the homotopy classes
 $$
 \xymatrix {
 v(p,q)\:
 S^{n-1}\times S^n
 \ar@{~>}[rr]^-{[b]} &&
 S^n\vee S^{2n-1}
 \ar@{~>}[rr]^-{(pi,qj)} &&
 S^n
 }
 $$
 (wavy arrows present homotopy classes)
 and
 $\=v(p,q)=[l]\circ v(p,q)\in[S^{n-1}\times S^n,S^n_\Q]$.
 Obviously,
 $v(0,q)=u(q)$ and
 $\=v(0,q)=\=u(q)$.
 We have
 $v(p,q)=v(p,0)$ if $p\divs q$
 (the proof is omitted) and
 $\=v(p,q)=\=v(p,0)$ if $p\ne0$
 \cite[Example~4.6]{ArkLup}.

 \begin {demo} [Proof of 1.5.]
 Take $q\in\Z$.
 Put $W=S^n\vee\dotso\vee S^n\vee S^{2n-1}$
 ($r$ times $S^n$).
 Let $d\:S^n\vee S^{2n-1}\to W$ be a map with
 $[d]=([\incl^W_1]+\dotso+[\incl^W_r],[\incl^W_{r+1}])$.
 Put $k=d\circ b\:S^{n-1}\times S^n\to W$.
 Let $h\:W\to S^n$ be a map with $[h]=(i,\dotsc,i,qj)$.
 By Lemma~12.1,
 $$
 \sum_{e\in\E^{r+1}}
 (-1)^{\*e}f([h\circ\M^W_e\circ k])=0.
 $$
 Since $[h\circ\M^W_e\circ k]=v(e_1+\dotso+e_r,e_{r+1}q)$,
 we have
 $$
 \sum_{e'\in\E^r}
 (-1)^{\*{e'}}
 \sum_{e''\in\E}
 (-1)^{e''}
 f(v(\*{e'},e''q))=0.
 $$
 Assume $r!\divs q$.
 If $e'\ne(0,\dotsc,0)$,
 the inner sum vanishes because
 then
 $\*{e'}\divs q$ and,
 consequently,
 the class $v(\*{e'},e''q)$ does not depend on $e''$.
 We get $f(v(0,0))-f(v(0,q))=0$,
 i.~e.\ $f(u(q))=f(u(0))$.
 \qed
 \end {demo}

 \begin {demo} [Proof of 1.6.]
 Assume $\Deg f\le r\in\N$.
 Take $q\in\Z$.
 As in the proof of 1.5,
 we get
 $$
 \sum_{e'\in\E^r}
 (-1)^{\*{e'}}
 \sum_{e''\in\E}
 (-1)^{e''}
 f(\=v(\*{e'},e''q))=0.
 $$
 If $e'\ne(0,\dotsc,0)$,
 the class $\=v(\*{e'},e''q)$ does not depend on $e''$.
 As in the proof of 1.5,
 we get $f(\=u(q))=f(\=u(0))$.
 \qed
 \end {demo}

 \begin {demo} [Proof of 1.7.]
 Assume $\Deg f\le r\in\N$.
 Consider the invariant $\~f\:\pi_{2n-1}(S^n)\to\Q$,
 $\~u\mapsto f(\~u\circ[c])$.
 By Lemma~10.1,
 $\Deg\~f\le r$.
 By Lemma~12.2,
 $\~f$ is gentle.
 Consider the function $F\:\Z\to\Q$, $q\mapsto f(u(q))$.
 We have $F(q)=\~f(qj)$.
 Therefore
 $F$ is gentle,
 i.~e.,
 by Lemma~3.11,
 is given by a polynomial.
 By 1.5,
 $F(q)=F(0)$ if $r!\divs q$.
 It follows that
 $F$ is constant.
 \qed
 \end {demo}



 {\noindent \tt ssp@pdmi.ras.ru}

 {\noindent \tt http://www.pdmi.ras.ru/\"{}ssp}

 \end {document}